\documentclass[12pt]{amsart}
\usepackage{amssymb}

\newtheorem{Thm}{Theorem}[section]

\newtheorem{corollary}[Thm]{Corollary}

\newtheorem{proposition}[Thm]{Proposition}
\newtheorem{definition}[Thm]{Definition}

\newtheorem{theorem}[Thm]{Theorem}

\newtheorem*{kadison-singer}{Kadison-Singer Problem}
\newtheorem*{paving conjecture}{Paving Conjecture}
\newtheorem*{bourgain-tzafriri}{Bourgain-Tzafriri Conjecture}
\newtheorem*{R-epsilon}{$R_{\epsilon}$-Conjecture}
\newtheorem*{FR-epsilon}{Finite $R_{\epsilon}$-Conjecture}

\def\ldots{\mathinner{\ldotp\ldotp\ldotp}}
\def\ldots{\mathinner{\cdotp\cdotp\cdotp}}

\def \D{{\mathbb D}}
\def \H{\mathbb H}
\def \N{\mathbb N}

\def \cal{\mathcal}

\def \beq{\begin{eqnarray*}}
\def \eeq{\end{eqnarray*}}

\begin{document}

\title[The Paving Conjecture and Triangular
matrices]{the paving conjecture is equivalent to the paving 
conjecture for triangular matrices}
\author[P.G. Casazza and J.C. Tremain
 ]{Peter G. Casazza and Janet C. Tremain}
\address{Department of Mathematics, University
of Missouri, Columbia, MO 65211-4100}

\thanks{The first author was supported by NSF DMS 0405376}

\email{janet,pete@math.missouri.edu}




\maketitle

\begin{abstract}
We resolve a 25 year old problem by showing that
  The Paving Conjecture is equivalent to
The Paving Conjecture for Triangular Matrices.
\end{abstract}

\section{Introduction}\label{Intro}
\setcounter{equation}{0}

The Kadison-Singer Problem \cite{KS} has been one of the most intractable
problems in mathematics for nearly 50 years.  

\begin{kadison-singer}[KS]
Does every pure state on the (abelian) von Neumann algebra $\D$ of
bounded diagonal operators on ${\ell}_2$ have a unique extension to
a (pure) state on $B({\ell}_2)$, the von Neumann algebra of all
bounded linear operators on the Hilbert space ${\ell}_2$?
\end{kadison-singer}

A {\bf state} of a von Neumann algebra ${\cal R}$ is a linear
functional $f$ on ${\cal R}$ for which $f(I) = 1$ and $f(T)\ge 0$
whenever $T\ge 0$ (i.e. whenever $T$ is a positive operator).
The set of states of ${\cal R}$ is a convex subset of the dual space
of ${\cal R}$ which is compact in the $w^{*}$-topology.  By the
Krein-Milman theorem, this convex set is the closed convex hull of its
extreme points.  The extremal elements in the space of states are
called the {\bf pure states} (of ${\cal R}$).  The Kadison-Singer
Problem had been dorment for many years when it was recently 
 brought back to life in \cite{CT} and \cite{CFTW} where
it was shown that
KS is equivalent to fundamental unsolved problems in a dozen different
areas of research in pure mathematics, applied mathematics and 
engineering.

A significant advance on the Kadison-Singer Problem
 was made by Anderson \cite{A3} in
1979 when he reformulated KS into what is now known as the
{\bf Paving Conjecture} (Lemma 5 of \cite{KS} shows a connection
between KS and Paving).  Before we state this conjecture, let us
introduce some notation.  For an operator $T$ on $\ell_2^n$,
 its matrix representation
$(\langle Te_i,e_j\rangle )_{i,j\in I}$ 
is with
respect to the natural orthonormal basis.
If $A\subset \{1,2,\ldots ,n\}$, 
the {\bf diagonal projection}
$Q_A$ is the matrix all of whose entries are zero except
for the $(i,i)$ entries for $i\in A$ which are all one.

\begin{paving conjecture}[PC]
For $\epsilon >0$, there is a natural number $r$ so
that for every
 natural number $n$ and
every linear operator
  $T$ on $l_2^n$
whose matrix has zero diagonal, 
  we can find a partition (i.e. a {\it paving})
$\{{A}_j\}_{j=1}^r$
  of $\{1, \ldots, n\}$, so that
  $$
  \|Q_{{A}_j} T Q_{{A}_j}\|  \le  \epsilon \|T\|
  \ \ \ \text{for all $j=1,2,\ldots ,r$.}
  $$
\end{paving conjecture}

It is important that $r$ not depend on $n$ in PC.  
We will say that an arbitrary operator $T$ satisfies PC if
$T-D(T)$ satisfies PC where $D(T)$ is the diagonal of $T$.
It is known that the class of
operators satisfying PC (the {\bf pavable operators})
is a closed subspace of $B({\ell}_2)$.  Also,
to verify PC we only need to verify it for
any one the following classes of operators \cite{AA,CFTW,CEKP}:
1.  unitary operators,
2.  positive operators,  3.  orthogonal projections (or just orthogonal
projections with $1/2's$ on the diagonal),
4.  Gram operators of the form
$T^{*}T=(\langle f_i,f_j\rangle )_{i,j\in I}$ where $\|f_i\|=1$ 
and $Te_i = f_i$ is a bounded operator.  The only large
classes of operators which have been shown to be pavable are
``diagonally dominant'' matrices \cite{BCHL,BCHL2,G},
 matrices with all entries
real and positive \cite{BHKW,HKW} and matrices with small entries
\cite{BT}.

Since the beginnings of the {\it paving era}, it has been a natural
question whether PC is equivalent to PC for triangular operators 
This question was formally asked several times at meetings by Gary Weiss and
Lior Tzafriri and
appeared (for a short time) on the AIM website (http://www.aimath.org/The 
Kadison-Singer Problem) as an important question for PC.
In this paper we will verify this conjecture.  Given two conjectures
$C_1,\ C_2$ we say that $C_1$ {\bf implies} $C_2$ if a positive answer
to $C_1$ implies a positive answer for $C_2$.  They are {\bf equivalent}
if they imply each other. 

\section{Preliminaries}\label{Prelim}
\setcounter{equation}{0}
Recall that a family of
vectors $\{f_i\}_{i\in I}$ is a {\bf Riesz basic sequence} in
a Hilbert space $\H$ if there are constants $A,B>0$ so that
for all scalars $\{a_i\}_{i\in I}$ we have:
$$
A^2\sum_{i\in I}|a_i |^2 \le \|\sum_{i\in I}a_i f_i \|^2 \le
B^2\sum_{i\in I}|a_i|^2.
$$
We call $A,B$ the {\bf lower and upper Riesz basis
bounds} for $\{f_i\}_{i\in I}$.  If $\epsilon >0$ and
$A = 1-\epsilon, B=1+\epsilon$ we call $\{f_i\}_{i\in I}$ an
$\epsilon$-{\bf Riesz basic sequence}.  If $\|f_i\|=1$ for all
$i\in I$ this is a {\bf unit norm} Riesz basic sequence.
A natural question is whether we can improve the Riesz basis
bounds for a unit norm Riesz basic sequence by partitioning the sequence
into subsets.

\begin{R-epsilon}
For every $\epsilon >0$, every unit norm Riesz basic sequence
is a finite union of $\epsilon$-Riesz basic sequences. 
\end{R-epsilon}

The $R_{\epsilon}$-Conjecture was posed by Casazza and Vershynin \cite{CV}
where it was shown that KS implies this conjecture.
It is now known that the $R_{\epsilon}$-Conjecture is
equivalent to KS \cite{CT}.
We will show that PC for triangular operators implies a positive
solution to the $R_{\epsilon}$-Conjecture.  Actually, we need the finite 
dimensional
quantative version of this conjecture.   

\begin{FR-epsilon}
Given $0< \epsilon, A,B$, there is a natural number $r=r(\epsilon,
A,B)$ so that for every $n\in \N$ and every unit norm Riesz basic
sequence $\{f_i\}_{i=1}^{n}$ for $\ell_2^n$ with Riesz basis bounds
$0<A\le B$, there is a partition $\{A_j\}_{j=1}^{r}$ of $\{1,2,\ldots ,n\}$
so that for all $j=1,2,\ldots ,r$ the family $\{f_i\}_{i\in A_j}$ is 
an $\epsilon$-Riesz basic sequence.
\end{FR-epsilon}

 There are standard methods for
turning infinite dimensional results into quantative finite dimensional
results so we will just outline the proof of their equivalence.  We will
need a proposition from \cite{CCLV}.

\begin{proposition}\label{Prop1}
Fix a natural number $r$ and
assume for every natural number $n$ we have a partition
$\{A_{i}^{n}\}_{i=1}^{r}$ of $\{1,2,\ldots , n\}$.  Then
there are natural numbers $\{n_{1}<n_{2}<\cdots\}$ so that
if $j\in A_{i}^{n_{j}}$ for some $i \in \{ 1, \ldots, r\}$, 
then $j\in A_{i}^{n_{k}}$, for all
$k\ge j$.  Hence, if $A_{i} = \{j\ |\ j\in A_{i}^{n_{j}}\}$ then

(1)  $\{A_{i}\}_{i=1}^{r}$ is a partition of $\mathbb N$.

(2)  If $A_{i} = \{j_{1}^i<j_{2}^i<\cdots \}$ then for every natural
number $k$ we have $\{j_{1}^i,j_{2}^i,\ldots , j_{k}^i\}\subset
A_{i}^{n_{j_{k}}}$.
\end{proposition}

\begin{theorem}
The $R_{\epsilon}$-Conjecture is equivalent to the Finite 
$R_{\epsilon}$-Conjecture.
\end{theorem}

\begin{proof}
Assume the Finite $R_{\epsilon}$-Conjecture is true.  Let 
$\{f_i\}_{i=1}^{\infty}$
be a unit norm 
Riesz basic sequence in $\H$ with bounds $0<A,B$ and fix $\epsilon >0$.
Then there is a natural number $r\in \N$ so that for all $n\in \N$ 
there is a partition $\{A_j^n\}_{j=1}^{r}$ of $\{1,2,\ldots ,n\}$ and
for every $j=1,2,\ldots ,r$ the
family $\{f_i\}_{i\in A_j^n}$ is an $\epsilon$-Riesz
basic sequence.  Choose a partition 
$\{A_j\}_{j=1}^r$ of $\N$ satisfying
 Proposition \ref{Prop1}.  By (2) of this proposition, 
for each $j=1,2,\ldots r$, the first $n$-elements of $\{f_i\}_{i\in A_j}$ 
come from one of the $A_{\ell}^m$ and hence form an $\epsilon$-Riesz
basic sequence.  So $\{f_i\}_{i\in A_j}$ is an $\epsilon$-Riesz basic
sequence.

Now assume the the Finite $R_{\epsilon}$-Conjecture fails.  Then there
is some $0<\epsilon,A,B$, natural numbers $n_1<n_2< \cdots$ and unit
norm Riesz basic sequences $\{f_i^r\}_{i=1}^{n_r}$ for $\ell_2^{n_r}$ so 
that whenever $\{A_j\}_{j=1}^{r}$ is a partition of $\{1,2,\ldots ,n_r\}$
one of the sets $\{f_i^r\}_{i\in A_j}$ is not an $\epsilon$-Riesz
basic sequence.  Considering 
$$
\{f_i\}_{i=1}^{\infty} = \{f_i^r\}_{i=1,r=1}^{\ n_r , \ \infty}
\in \left ( \sum_{r=1}^{\infty}\oplus \ell_2^{n_r} \right )^{1/2},
$$
we see that this family of vectors forms a unit norm Riesz basic sequence
with bounds $0<A,B$ but for any natural number $r$ and any partition
$\{A_j\}_{j=1}^{r}$ of $\N$ one of the sets $\{f_i\}_{i\in A_j}$
is not an $\epsilon$-Riesz basic sequence.
\end{proof}

\section{The Main Theorem}\label{Main}
\setcounter{equation}{0}

Our main theorem is:

\begin{theorem}\label{T}
The Paving Conjecture is equivalent to the Paving Conjecture for Triangular
matrices.
\end{theorem}

\begin{proof}
Since a paving of $T$ is also a paving of $T^{*}$, we only need to show
that The Paving Conjecture for Lower Triangular Operators implies the
Finite $R_{\epsilon}$-Conjecture.  Fix $0< \epsilon, A,B$, fix
$n\in \N$ and let $\{f_i\}_{i=1}^{n}$ be a 
unit norm Riesz basis for $\ell_2^n$
with bounds $A,B$.  We choose a natural number $r\in \N$ satisfying:
$$
1-\frac{B^4}{A^4r} \ge 1- \frac{\epsilon}{2}.
$$
We will do the proof in 5 steps.
\vskip12pt
\noindent {\bf Step 1}:  There is a partition $\{A_j\}_{j=1}^{r}$ of 
$\{1,2,\ldots ,n\}$ so that for every $j=1,2,\ldots ,r$ and every
$i\in A_j$ and every $1\le k\not= j \le r$ we have:
$$
\sum_{i\not= \ell \in A_j}|\langle f_i,f_{\ell}\rangle |^2 \le
\sum_{\ell \in A_k}|\langle f_i,f_{\ell}\rangle |^2.
$$
\vskip12pt
The argument for this is due to Halpern, Kaftal and Weiss (\cite{HKW},
Proposition 3.1) so we will outline it for our case.
Out of all ways of partitioning $\{1,2,\ldots ,n\}$ into $r$-sets,
choose one, say $\{A_j\}_{j=1}^{r}$, which minimizes
\begin{equation}\label{E1}
\sum_{j=1}^{r}\sum_{i\in A_j}\sum_{i\not= \ell \in A_{j}}
|\langle f_i,f_{\ell}\rangle |^2.
\end{equation}
We now observe that for each $1\le j\le r$, each $i\in A_j$ and all 
$1\le k\not= j \le r$ we have
$$
\sum_{i\not= \ell \in A_j}|\langle f_i,f_{\ell}\rangle |^2 
\le \sum_{\ell \in A_k}|\langle f_i,f_{\ell}\rangle |^2.
$$
To see this, assume this inequality fails.  That is, for some $j_0,
i_0,k_0$ as above we have
$$
\sum_{i_{0}\not= \ell \in A_{j_0}}|\langle f_{i_0},f_{\ell}\rangle |^2 > 
\sum_{\ell \in A_{k_0}}|\langle f_{i_0},f_{\ell}\rangle |^2.
$$
We define a new partition $\{B_j\}_{j=1}^{r}$ of $\{1,2,\ldots , n\}$
by: $B_j=A_j$ if $j\not= j_0,k_0$; $B_{j_0} = A_{j_0}-\{i_0\}$;
$B_{k_0} = A_{k_0} \cup \{i_0\}$.  It easily follows that
$$
\sum_{j=1}^{r}\sum_{i\in B_j}\sum_{i\not= \ell \in B_j}
|\langle f_i,f_{\ell}\rangle |^2 < \sum_{j=1}^{r} \sum_{i\in A_j}
\sum_{i\not= \ell \in A_j}|\langle f_i,f_{\ell}\rangle |^2,
$$
which contradicts the minimality of Equation \ref{E1}.
\vskip12pt

\noindent {\bf Step 2}:  For every $j=1,2,\ldots ,r$ and every
$i\in A_j$ we have
$$
\sum_{i\not= \ell \in  A_j}|\langle f_i,f_{\ell}\rangle |^2 \le 
\frac{B^2}{r}.
$$
\vskip12pt
Define an operator $Sf = \sum_{i=1}^{n} \langle f,f_i\rangle f_i$.  Then,
$$
\langle Sf,f\rangle = \sum_{i=1}^{r} |\langle f,f_i\rangle |^2,
$$
and since $\{f_i\}_{i=1}^{n}$ is a Riesz basis with bounds $A,B$ we have
$$
A^2I \le S \le B^2I.
$$
Now, by Step 1,
\begin{eqnarray*}
\sum_{i\not= \ell \in A_j}|\langle f_i,f_{\ell}\rangle |^2 &\le&
\frac{1}{r} \left [ \sum_{i\not= \ell \in A_j}|\langle f_i,f_{\ell}\rangle
|^2 + \sum_{j\not= k=1}^{r} \sum_{\ell \in A_k}|\langle f_i,f_{\ell}
\rangle |^2 \right ]\\
&\le& \frac{1}{r}\sum_{i=1}^{n}|\langle f_i,f_{\ell}\rangle |^2\\
&\le& \frac{1}{r}\|S\|\|f_i\|^2 \le \frac{B^2}{r}.
\end{eqnarray*}
\vskip12pt
\noindent {\bf Step 3}:  For each $j=1,2,\ldots ,r$ and all $i\in A_j$,
if $P_{ij}$ is the orthogonal projection of span $\{f_{\ell}\}_{\ell
\in A_j}$ onto 
span $\{f_{\ell}\}_{i\not= \ell \in A_j}$ then
$$
\|P_{ij}f_i\|^2 \le \frac{B^4}{A^4r}.
$$
\vskip12pt
Define the operator $S_{ij}$ on span $\{f_{\ell}\}_{\ell \in A_j}$ by
$$
S_{ij}(f) = \sum_{i\not= \ell \in A_j}\langle f,f_{\ell}\rangle f_{\ell}.
$$
Then $A^2I \le S_{ij}\le B^2I$ and 
$\{S_{ij}^{-1}f_{\ell}\}_{i\not= \ell \in A_j}$ are the dual
functionals for the Riesz basic sequence $\{f_{\ell}\}_{i\not= \ell
\in A_j}$.  Also, as in Step 1, $A^2I \le S_{ij}\le B^2I$.  So by Step 2,
\begin{eqnarray*}
\|P_{ij}f_i\|^2 &=& \|\sum_{i\not= \ell \in A_j}\langle f_i,f_{\ell}
\rangle S_{ij}^{-1}f_{\ell}\|^2\\
&\le&\|S_{ij}^{-1}\|^2 \|\sum_{i\not= \ell \in A_j}
\langle f_i,f_{\ell}\rangle f_{\ell}\|^2\\
&\le& \frac{B^2}{A^4}\sum_{i\not= \ell \in A_j}|\langle f_i,f_{\ell}
\rangle |^2| \le \frac{B^4}{A^4 r}.
\end{eqnarray*}
\vskip12pt
\noindent {\bf Step 4}:  Fix $1\le j\le r$ and let $A_j = \{i_1,i_2,
\ldots ,i_k\}$.  If we Gram-Schmidt $\{f_{i_{\ell}}\}_{\ell =1}^{k}$
to produce an orthonormal basis $\{e_{i_{\ell}}\}_{\ell =1}^{k}$
then for all $1\le m\le k$ we have
$$
|\langle f_{i_m},e_{i_m}\rangle |^2 \ge 1- \frac{\epsilon}{2}.
$$
\vskip12pt
Fix $1\le m\le k$ and let $Q_m$ be the orthogonal projection of
span $\{e_{i_{\ell}}\}_{\ell = 1}^{k}$ onto 
span $\{e_{i_{\ell}}\}_{\ell =1}^{m}$ = span $\{f_{i_{\ell}}\}_{\ell =1}^{m}$.
By Step 3,
$$
\|Q_{m}f_{i_m}\|^2 \le \|P_{mj}f_{i_m}\|^2 \le \frac{B^4}{A^4r}.
$$
Since 
$$
f_{i_m} = \sum_{\ell =1}^{m} \langle f_{i_{\ell}},e_{i_{\ell}}\rangle e_{i_{\ell}},
$$
we have
\begin{eqnarray*}
|\langle f_{i_m},e_{i_m}\rangle |^2 &=& \|f_{i_m}\|^2 - 
\|Q_{m-1}f_{i_m}\|^2\\
&\ge& 1-\frac{B^4}{A^4r} \ge 1-\frac{\epsilon}{2},
\end{eqnarray*}
where the last inequality follows from our choice of $r$.
\vskip12pt

\noindent {\bf Step 5}:  We complete the proof.
\vskip12pt
Let
$$
M = \left ( \langle f_{i_s},e_{i_t}\rangle \right ) _{s\not= t =1}^{k},
$$
where by this notation we mean the $k\times k$-matrix with zero
diagonal and the given values off the diagonal.  By the Gram-Schmidt
Process, $M$ is a lower triangular matrix with zero diagonal.
Define an operator 
$T:\ell_2^k \rightarrow span \ \{e_{i_{\ell}}\}_{\ell =1}^{k}$ by
$$
T\left ( (a_{i_{\ell}})_{\ell = 1}^{k} \right ) = 
\sum_{\ell =1}^{k} a_{i_{\ell}}f_{i_{\ell}}.
$$
If $K$ is the matrix of $T$ with respect to the orthonormal
basis $\{e_{i_{\ell}}\}_{\ell =1}^{k}$ and $D=D(K)$ is the
diagonal of $K$ then $M = (K-D)^{*}$ and so 
$$
\|M\| \le \|K\| + \|D\| = \|T\|+1 \le B+1.
$$
By The Paving Conjecture for lower triangular matrices, there is
a natural number $L_j$ (which is a function of $0 < \epsilon$ and
$B$ only) and a partition $\{B_{\ell}^{j}\}_{\ell =1}^{L_j}$ of 
$\{i_1,i_2 \ldots , i_k\}$ so that
$$
\|Q_{B_{\ell}^j}MQ_{B_{\ell}^j}\| \le \frac{\epsilon}{2},
$$
for all $\ell = 1,2, \ldots ,L_j$ ($Q_{B_{\ell}^j}$ was defined in the
introduction).  Now, for all scalars $(a_{i_s})_{i_s \in B_{\ell}^j}$,
if 
$$
f = \sum_{i_s \in B_{\ell}^j}a_{i_s}f_{i_s}, 
$$
then
\begin{eqnarray*}
\|\sum_{i_s \in B_{\ell}^j}a_{i_s}f_{i_s}\|
&=& \| D(f) + Q_{B_{\ell}^j}M^{*}Q_{B_{\ell}^j}(f)
\|\\
&\ge& \|Df\| - \|Q_{B_{\ell}^j}M^{*}Q_{B_{\ell}^j}(f)\|\\
&\ge& (1-\frac{\epsilon}{2})\|f\| - \frac{\epsilon}{2}\|f\|\\
&\ge& (1-\epsilon)\|f\|.
\end{eqnarray*}
Similarly, 
$$
\|\sum_{i_s \in B_{\ell}^j}a_{i_s}f_{i_s}   \| \le (1+\epsilon)\|f\|.
$$  
It follows that $\{f_i\}_{i\in B_{\ell}^j}$ is an $\epsilon$-Riesz
basic sequence for all $j=1,2,\ldots,r$ and all $\ell = 1,2,\ldots ,L_j$. 
Hence, the Finite $R_{\epsilon}$-Conjecture holds which 
completes the proof of the theorem.
\end{proof}

Let us make an observation concerning the proof of the main theorem.

\begin{definition}
Let $\{f_i\}_{i=1}^{\infty}$ be a sequence of vectors
in a Hilbert space $\H$.  For each $i=1,2,\ldots $ let
$P_i$ be the orthogonal projection of $\H$ onto span 
$\{f_{\ell}\}_{i\not= \ell \in \N}$.  Our sequence is said to
be {\bf $\epsilon$-minimal} if $\|P_i\| \le \epsilon$ for all $i=1,2,\ldots$.
\end{definition}

The first three steps of the proof of Theorem \ref{T} yields:

\begin{corollary}
If $\{f_i\}_{i=1}^{\infty}$ is a unit norm Riesz basic sequence in a Hilbert
space $\H$ then for every $\epsilon >0$ there is a partition
$\{A_j\}_{j=1}^{r}$ of $\N$ so that for all $j=1,2,\ldots ,r$,
the family $\{f_i\}_{i\in A_j}$ is $\epsilon$-minimal.
\end{corollary}


\begin{thebibliography}{10}
\bibitem{AA}  C.A.  Akemann and J. Anderson, {\em Lyapunov theorems
for operator algebras}, Mem. AMS {\bf 94} (1991).

\bibitem{A3}  J. Anderson, {\em Extreme points in sets of positive
linear maps on $B(\H)$}, Jour. Functional Analysis {\bf 31} (1979)
195--217.

\bibitem{BCHL}  R. Balan, P.G. Casazza, C. Heil and Z. Landau,
{\em Density, overcompleteness and localization of frames.
I. Theory}, Preprint.

\bibitem{BCHL2}  R. Balan, P.G. Casazza, C. Heil and Z. Landau,
{\em Density, overcompleteness and localization of frames.
II. Gabor systems}, Preprint.

\bibitem{BHKW} K. Berman, H. Halpern, V. Kaftal and G. Weiss,
{\em Matrix norm inequalities and the relative Dixmier property},
Integ. Eqns. and Operator Theory {\bf 11} (1988) 28--48.


\bibitem{BT}  J. Bourgain and L. Tzafriri, {\em Invertibility
of ``large'' submatrices and applications to the geometry of
Banach spaces and Harmonic Analysis}, Israel J. Math.
{\bf 57} (1987) 137--224.

\bibitem{CCLV}  P.G. Casazza, O. Christensen, A. Lindner and
R. Vershynin, {\em Frames and the Feichtinger conjecture},
Proceedings of AMS, {\bf 133} No. 4 (2005) 1025--1033.

\bibitem{CEKP}  P.G. Casazza, D. Edidin, D. Kalra and V. Paulsen,
The Kadison-Singer Problem and Projections, Preprint.

\bibitem{CT} P.G. Casazza and J.C. Tremain, {\em The Kadison-Singer
Problem in Mathematics and Engineering}, Proceedings of the National
Academy of Sciences, {\bf 103} No. 7 (2006) 2032-2039.

\bibitem{CFTW}  P.G. Casazza, M. Fickus, J.C. Tremain, and
E. Weber, {\em The Kadison-Singer Problem in Mathematics
and Engineering:  Part II:  A detailed account}. 
(Accepted for The Proceedings of The 2005 Great Plains Operator Theory
Symposium (GPOTS), Contemp. Math., Amer. Math. Soc., to appear in 2006).

\bibitem{CV}  P.G. Casazza and R. Vershynin, {\em Kadison-Singer
meets Bourgain-Tzafriri}, Preprint.

\bibitem{G}  K.H. Gr\"ochenig, {\em Localized frames are finite unions
of Riesz sequences}, Adv. Comp. Math. {\bf 18} (2003) 149--157.

\bibitem{HKW} H. Halpern, V. Kaftal and G. Weiss,
{\em Matrix pavings and Laurent operators}, J. Op. Th. {\bf 16}
 (1986) 121--140.

\bibitem{KS} R. Kadison and I. Singer,
{\em Extensions of pure states}, American Jour. Math. {\bf 81}
(1959), 383--400.

\end{thebibliography}
\end{document}